\documentclass[11pt]{amsart}
\usepackage{amsfonts,amssymb,amsthm,amsmath,amsxtra,amscd,verbatim,eucal,latexsym}
\usepackage[all]{xy}
\usepackage[dvips]{graphics}

\theoremstyle{plain}
\newtheorem{thm}{Theorem}[section]
\newtheorem{lem}[thm]{Lemma}
\newtheorem{prop}[thm]{Proposition}
\newtheorem{cor}[thm]{Corollary}

\theoremstyle{definition}
\newtheorem{defi}[thm]{Definition}

\theoremstyle{remark}

\newtheorem{rmk}[thm]{Remark}

\newtheorem*{ack}{Acknowledgements}

\numberwithin{equation}{section}

\def\Z{\mathbb Z}
\def\N{{\mathbb N}}

\def\C{{\mathbb C}}
\def\A{{\mathbb A}}

\def\Q{{\mathbb Q}}
\def\P{{\mathbb P}}
\def\L{{\mathbb L}}

\def\CC{\mathcal{C}}

\def\O{{\mathcal O}}

\def\MM{\mathcal{M}}
\def\NN{\mathcal{N}}

\def\SS{\mathcal{S}}

\def\UU{\mathcal{U}}
\def\XX{\mathcal{X}}
\def\YY{\mathcal{Y}}

\def\1{{\bf 1}}

\def\a{\alpha}

\def\g{\gamma}

\def\ff{\psi}
\def\e{\eta}

\def\n{\nu}
\def\m{\mu}
\def\om{\omega}
\def\p{\pi}
\def\r{\rho}
\def\s{\sigma}
\def\t{\tau}
\def\q{\chi}
\def\x{\xi}
\def\z{\zeta}
\def\D{\Delta}
\def\G{\Gamma}
\def\LLL{\Lambda}

\def\.{\cdot}
\def\^{\widehat}
\def\~{\widetilde}
\def\o{\circ}
\def\ov{\overline}

\def\({\left(}
\def\){\right)}

\DeclareMathOperator{\SL}  {SL}

\DeclareMathOperator{\Spec} {Spec}
\DeclareMathOperator{\Aut} {Aut}
\DeclareMathOperator{\Space} {\mathfrak{Space}}
\DeclareMathOperator{\Sch} {\mathfrak{Sch}}

\DeclareMathOperator{\DM} {\mathfrak{DM}}

\DeclareMathOperator{\Cont} {Cont}
\DeclareMathOperator{\Ex} {Ex}
\DeclareMathOperator{\ord} {ord}
\DeclareMathOperator{\SM} {SM}
\DeclareMathOperator{\st} {str}

\DeclareMathOperator{\Cyl} {Cyl}
\DeclareMathOperator{\age} {age}

\newcommand{\coarse}[1]{\ensuremath{{ \overline{#1}}}}

\begin{document}

\title{Stringy Chern classes of singular varieties}

\thanks{The first author was partially supported by
the University of Michigan Rackham Research Grant and the MIUR of
the Italian Government, National Research
Project ``Geometry on Algebraic Varieties" (Cofin 2002).
The second author was partially supported by CONACYT.
The third author was partially supported by an NSF Postdoctoral
Fellowship at the University of Michigan.
Research of the fourth author was partially carried out during
his visit at the Max Planck Institut, Bonn.}

\author[T. de Fernex]{Tommaso de Fernex}
\address{Department of Mathematics, University of Michigan, Ann Arbor,
MI 48109, USA}
\email{defernex@umich.edu}
\curraddr{Institute for Advanced Study, School of Mathematics,
1 Einstein Drive, Princeton, NJ 08540, USA}

\author[E. Lupercio]{Ernesto Lupercio}
\address{Departamento de Matem\'aticas\\CINVESTAV\\C.P. 07000, M\'exico, D.F.,
M\'exico} \email{lupercio@math.cinvestav.mx}

\author[T. Nevins]{Thomas Nevins}
\address{Department of Mathematics, University of Illinois, Urbana,
IL 61801, USA}
\email{nevins@math.uiuc.edu}

\author[B. Uribe]{Bernardo Uribe}
\address{Departamento de Matem\'aticas, Universidad de los Andes, Bogot\'a,
Colombia} \email{buribe@uniandes.edu.co}

\begin{abstract}
Motivic integration \cite{Kon} and MacPherson's transformation \cite{Mac}
are combined in this paper to construct a theory of ``stringy" Chern classes
for singular varieties. These classes enjoy strong birational invariance properties,
and their definition encodes data coming from resolution of singularities.
The singularities allowed in the theory
are those typical of the minimal model program;
examples are given by quotients of manifolds by finite groups.
For the latter an explicit formula is proven,
assuming that the canonical line bundle of the manifold descends to the quotient.
This gives an expression of the stringy Chern class of the quotient
in terms of Chern-Schwartz-MacPherson classes
of the fixed-point set data.
\end{abstract}

\maketitle

\section*{Introduction}

Chern classes provide a powerful tool in the study of manifolds,
and the problem of finding good notions of Chern classes for
singular varieties has been the object of works by several
authors, including \cite{Sch,Mac,Ful,FJ,Tot,Alu}.

The purpose of this paper is to introduce and study
a new theory of Chern classes for singular complex varieties by
combining techniques of motivic integration
with the transformation defined by MacPherson in \cite{Mac}.
In particular, we give a construction of a class, which we call the
{\it stringy Chern class},
with strong birational properties and an interesting
``orbifold'' interpretation when
the variety is a Gorenstein quotient of a manifold by a finite group.
The adjective ``stringy" here
is dictated by the analogy with other invariants that have been
introduced in the last decade or so,
see for example \cite{DHVW,Zas,Bat1,Bat2,BL2}.

The notion of {\it $K$-equivalence} of varieties plays an important role in
birational geometry \cite{Kaw} (see also Definition~\ref{defi-K} below).
It was proven by Aluffi \cite{Alu} that, after passing to rational
coefficients, the
Chern classes of two $K$-equivalent smooth varieties are the
push-forward of the same class on a common resolution.
Consequently, given a possibly
singular variety $X$ (with suitable restrictions on
the singularities), one would like to define a class $c_{\st}(X) \in A_*(X)_\Q$
satisfying:
\begin{enumerate}
\item[(i)]
$c_{\st}(X) = c(T_X) \cap [X]$ if $X$ is smooth, and
\item[(ii)]
if $X$ and $X'$ are $K$-equivalent and
$$
\xymatrix{
& Y \ar[dl]_f \ar[dr]^{f'} \\
X && X'
}
$$
is a common resolution, then $c_{\st(X)} = f_*C$ and $c_{\st(X')} = f'_*C$
for some $C \in A_*(Y)_\Q$.
\end{enumerate}

Let us briefly recall the definition of the Chern-Schwartz-MacPherson class of
\cite{Sch, Mac}, following \cite{Mac}. Given an arbitrary complex variety $X$, MacPherson
constructs a group homomorphism
$$
c : F(X) \to A_*(X)
$$
from the group of constructible functions on $X$
to the Chow group of $X$. This transformation
is functorial with respect to proper morphisms
and satisfies $c(\1_X) = c(T_X) \cap [X]$ whenever $X$ is smooth.
Then the {\it Chern-Schwartz-MacPherson class} of a complex variety $X$
is the class $c_{\SM}(X) := c(\1_X) \in A_*(X)$.
While the class $c_{\SM}(X)$ satisfies property~(i) above,
simple examples (such as a surface with rational double points, compared
to its minimal resolution) show that it does not enjoy property~(ii).

In this paper, we encode information coming from resolution
of singularities into a (${\mathbb Q}$-valued)
constructible function $\Phi_X$ that
we use in place of $\1_X$.
The {\it stringy Chern class} of $X$
is then defined to be the image of this function
via the group homomorphism $c : F(X)_\Q \to A_*(X)_\Q$
obtained from MacPherson's transformation, that is, the class
$$
c_{\st}(X) := c(\Phi_X) \in A_*(X)_\Q.
$$
\begin{thm}
The stringy Chern class satisfies properties~(i) and~(ii).
\end{thm}

Our construction of the function $\Phi_X$
uses the powerful and flexible technique of
motivic integration
\cite{Kon, Bat2, DL1, Loo, Rei}.
Motivic integration has by now been employed to extend several
numerical invariants of manifolds to singular varieties.
By combining the MacPherson
transformation with motivic integration (with its natural change-of-variables
formula), we obtain a direct
construction of the cycle class $c(\Phi_X)$
with manifest birational invariance properties.

More precisely, our construction of the function $\Phi_X$ takes the following
form. Given a complex variety $X$, we define a natural
ring homomorphism $\Phi$ from suitable relative motivic rings over $X$
to the ring $F(X)_\Q$ of rational-valued constructible functions on $X$.
Then we restrict our attention to the case of a normal variety $X$ with a
$\Q$-Cartier canonical divisor and having at most log-terminal
singularities.\footnote{These are, for instance, the singularities
of the quotient of a smooth variety by the action of a finite group;
in general they form a natural class of singularities in birational
geometry.}
We fix a resolution of singularities $f : Y \to X$, and denote
by $K_{Y/X}$ the relative canonical $\Q$-divisor. Then we define
$$
\Phi_X := \Phi\left(\int_{Y_\infty}\L_X^{-\ord(K_{Y/X})} d\mu^X\right)
\in F(X)_\Q,
$$
where the argument of $\Phi$ in the right hand side is a
motivic integral computed relative to $X$.
Standard properties of motivic integration ensure that
this function is independent of the resolution and that the
stringy Chern class
$$
c_{\st}(X) := c(\Phi_X) \in A_*(X)_\Q
$$
satisfies properties~(i) and (ii), as desired.
Moreover,
whenever $\Phi_X$ is an integral-valued function, we can apply directly
the homomorphism $c : F(X) \to A_*(X)$ constructed by MacPherson
without passing to rational coefficients, thus defining a class
$$
\~c_{\st}(X) := c(\Phi_X) \in A_*(X),
$$
still satisfying (i) and (ii), whose image in $A_*(X)_\Q$ is equal to $c_{\st}(X)$.

The connection between stringy Chern classes and other stringy invariants
is manifested when $X$ is proper, as in this case
the degree of $c_{\st}(X)$ is equal to
the stringy Euler number $e_{\st}(X)$ defined by Batyrev \cite{Bat2}.
Furthermore, we establish a formula for quotient varieties,
where stringy Chern classes
are explicitly computed in terms of Chern classes of fixed-point set data.
Suppose that $X$ is the quotient of a smooth
quasi-projective variety $M$ under the action of a
finite group $G$. We assume that the canonical line bundle of $M$
descends to $X$. This is for instance a typical situation
considered in mirror symmetry; the reader may also think of the case
of a finite subgroup of $\SL_n(\C)$ acting on $\C^n$.
Fix a set $\CC(G)$ of representatives of
conjugacy classes in $G$. For any $g \in G$, denote by
$M^g \subseteq M$ the corresponding fixed-point set, and by
$C(g) \subseteq G$ the centralizer of $g$ in $G$.
Note that there is a morphism $\p_g : M^g/C(g) \to X$ commuting
with the quotient maps. Then we prove the following formula.

\begin{thm}
With the above notation, $\Phi_X \in F(X)$ and
\begin{equation}\label{intro-MK}
\~c_{\st}(X) = \sum_{g \in \CC(G)} (\p_g)_*c_{\SM}(M^g/C(g)) \quad\text{in $A_*(X)$}
\end{equation}
\end{thm}

If $X$ is proper, then as a corollary of this result we obtain
(under our assumptions) Batyrev's formula \cite{Bat2}
for the stringy Euler number of a quotient variety.

The formula in~\eqref{intro-MK} acquires particular interest in the
language of stacks. More specifically, it turns out that there is a deep
connection between the theory of stringy Chern classes and the
theory of characteristic classes for Deligne-Mumford stacks.
A hint of such a connection
is given by Corollary~\ref{MK-orb} below, where the correspondence
is stated on the level of constructible functions.

\medskip

We close this introduction with a few remarks on the
techniques used in this paper in comparison to related works.

Motivic integration was invented by Kontsevich \cite{Kon} to prove Batyrev's
conjecture on the invariance of Hodge numbers for birational Calabi-Yau manifolds
(and in fact for $K$-equivalent smooth varieties),
and then applied to an array of different problems.
The basic techniques of motivic integration used in this paper
can be found in a series of works, for instance \cite{Bat1,DL1,DL2,Loo,Rei},
and for more introductory treatments we recommend \cite{Bli,Cra}.
An introduction to stringy invariants in the context of
motivic integration can be found in \cite{Vey}.

An alternative approach to the same kind of problems is offered by the
weak factorization theorem \cite{AKMW}.
This is the approach followed by Aluffi in \cite{Alu}.
In fact, the same approach was previously followed by Borisov and Libgober
to prove similar results for elliptic genera and elliptic classes
of manifolds and, moreover, to generalize
these invariants to singular varieties and determining
their ``McKay correspondence" \cite{BL2,BL4}.
It was explained to us by Borisov and Libgober that one can
reconstruct the stringy Chern class
defined in this paper from the orbifold elliptic class
defined in \cite{BL4} by taking the coefficients of a certain
Laurent expansion in one of the two variables after
a limiting process in the other variable.
In particular, this relationship gives evidence
of a possible motivic interpretation of the elliptic genus,
which we are currently exploring.

When the writing of this paper was almost complete,
we learned of a new preprint by Aluffi \cite{Alu2},
where stringy Chern classes are independently defined.
The definition of stringy Chern class given by Aluffi follows
a different route than ours.
His technique is based on the use of a category
called a {\it modification system} that permits
to study intersection-theoretic invariants in the
birational class of a given variety. Motivic
integration is not used in Aluffi's approach; nevertheless,
there are clear analogies between basic properties
of the two techniques.
Our approach to the definition of stringy Chern classes
also gives, a posteriori, a concrete bridge between Aluffi's theory
of modification systems and motivic integration.
More recently, the teory of stringy Chern classes
has been integrated in a more general framework by
Brasselet, Sch\"urmann and Yokura in their preprint
\cite{BSY2}, which is a completely new and improved version
of \cite{BSY}.

\begin{ack}
The first named author is grateful for an early conversation with
Paolo Aluffi on the results contained in \cite{Alu}, which has
very much motivated this work. We would like to thank Lev
Borisov, Pierre Deligne, Lawrence Ein, William Fulton, Paul Horja, Dano Kim,
Andrew Kresch, Robert Lazarsfeld, Anatoly Libgober and Willem Veys
for several useful comments and discussions,
and the referees for their remarks, suggestions and corrections.
\end{ack}

\section{Motivic Integration in the Relative Setting}

In this section we review the main features
of motivic integration over a fixed complex algebraic variety
(an integral separated scheme of finite type over $\C$).
A good reference source for motivic integration in the relative
setting is \cite{Loo}.
Knowledge of the basic concepts in the theory of arc spaces
will be assumed. The reader new to arc spaces and motivic integration
can find nice introductions to the subject in \cite{Cra,Bli}.

We start with the construction of the relative motivic ring.
We fix a complex algebraic variety $X$ (an integral separated scheme
of finite type over $\C$).
Let $\Sch_X$ be the category of separated schemes of finite type
over $X$, or {\it $X$-schemes} for short.
Given an $X$-scheme $g : V \to X$, we denote
by $\{g:V \to X\}$, or simply by $\{V\}$ when $X$ and $g$ are clear
from the context, the corresponding
class modulo isomorphism over $X$.
Moreover, we set $\L_X := \{\A^1_X\}$.

Let $K_0(\Sch_X)$ denote the free $\Z$-module generated by the
isomorphism classes of $X$-schemes, modulo the relations
$\{V\} = \{V \setminus W\} + \{W\}$ whenever
$W$ is a closed subscheme of an $X$-scheme $V$,
and both $W$ and $V\setminus W$ are viewed as $X$-schemes
under the restriction of the morphism $V \to X$.
$K_0(\Sch_X)$ becomes a ring when the product is defined by setting
$\{V\} \. \{W\} = \{V \times_X W\}$ and extending it associatively.
This ring has zero $\{\emptyset\}$ and for identity $\{X\}$.

\begin{rmk}\label{rmk:sch-v-var}
For every $X$-scheme $V$, we have $\{V\} = \{V_{\rm red}\}$ in
$K_0(\Sch_X)$. Indeed, we can stratify $V$ into affine strata $V_i$ (over $X$),
and given embedding $V_i \subseteq \A^{n_i}_X$, we obtain
$\{(V_i)_{\rm red}\} = \{\A^{n_i}_X\} - \{U_i\} = \{V_i\}$,
if $U_i$ is the complement of $V_i$ in $\A^{n_i}_X$, and
the $(V_i)_{\rm red}$ glue back together to form
a stratification of $V_{\rm red}$.
In particular, $K_0(\Sch_X)$ is generated by isomorphism classes of
$X$-varieties.
\end{rmk}

We define $\MM_X := K_0(\Sch_X)[\L_X^{-1}]$.
We will use the symbol $\{V\}$ also to denote the class
of a $X$-scheme $V$ in $\MM_X$.
The dimension $\dim \a$ of an element $\a \in \MM_X$
is by definition the infimum of the set of integers $d$ for
which $\a$ can be written as a finite sum
$$
\a = \sum m_i\{V_i\}\L_X^{-b_i}
$$
with $m_i \in \Z$ and $\dim V_i - b_i \le d$
(here the dimension of $V_i$ is the one over $\Spec \C$).
Note that $\dim \{\emptyset\} = -\infty$.
The dimension function satisfies
$\dim(\alpha +\beta)\leq \max\{\dim(\alpha),\dim(\beta)\}$ and
$\dim(\alpha\cdot\beta)\leq \dim(\alpha)+\dim(\beta)$,
so we obtain a structure of
filtered ring on $\MM_X$ with the filtration
of $\MM_X$ given by dimension.
Completing with respect to the dimensional
filtration (for $d \to -\infty$),
we obtain the {\it relative motivic ring} $\^\MM_X$.
We will use $\{V\}$ also to denote the image of $V$
in $\^\MM_X$ under the natural map $\MM_X \to \^\MM_X$.
We will denote by
$$
\t = \t_X : K_0(\Sch_X) \to \^\MM_X
$$
the composition of the maps $K_0(\Sch_X) \to \MM_X$ and $\MM_X \to \^\MM_X$.

\begin{rmk}
It is not known whether either one of the two factors
of $\t$ is injective.
\end{rmk}

All main definitions and properties valid for the motivic integration
over $\Spec \C$ translate to the relative setting by simply remembering
the maps over $X$. For instance, consider a nonsingular $X$-variety
$f : Y \to X$. Let $Y_\infty$ be the space of arcs of $Y$, and
denote by $\Cyl(Y_\infty)$ the set of cylinders on $Y_\infty$.
Then the motivic pre-measure
$$
\m^X : \Cyl(Y_\infty) \to \^\MM_X
$$
is defined as follows.
For any cylinder $C \in \Cyl(Y_\infty)$, we choose
an integer $m$ such that $\p_m^{-1}(\p_m(C)) = C$
(here $\p_m : Y_\infty \to Y_m$ is the truncation map to the space of $m$th jets).
Then we put
$$
\m^X(C) := \{\p_m(C)\}\L_X^{-m\dim Y},
$$
where $\p_m(C)$ is viewed as a constructible set over $X$ under the
composite morphism $Y_m \to Y \to X$. A standard computation
shows that the definition does not depend on the choice of $m$.

For any effective divisor $D$ on $Y$,
we denote by $\ord(D) : Y_\infty \to \N \cup \{\infty\}$ the
order function along $D$, and set
$\Cont^p(D) := \{\g \in Y_\infty \mid \ord_\g(D) = p\}$.
This is a cylinder in $Y_\infty$.
Then the relative motivic integral is defined by
\begin{equation}\label{int}
\int_{Y_\infty} \L_X^{-\ord (D)} d\m^X
:= \sum_{p \ge 0} \m^X(\Cont^p(D))\L_X^{-p}.
\end{equation}
This gives an element in $\^\MM_X$.

The following formula is a basic (but extremely useful)
property of motivic integration. A proof for integration
over $\Spec \C$ can be found in \cite{DL1}, and
the same proof translates in the relative setting
by keeping track of the morphisms to $X$ and observing
that $Y \times \A^1 = Y \times_X \A^1_X$ for any
$X$-variety $Y$.

\begin{thm}[Change of variables formula \cite{Kon}]\label{change-of-vars}
Let $g : Y' \to Y$ be a proper birational map between nonsingular
varieties over $X$, and let $K_{Y'/Y}$
be the relative canonical divisor of $g$.
Let $D$ be an effective divisor on $Y$.
Then
$$
\int_{Y_\infty} \L_X^{-\ord (D)} d\m^X =
\int_{Y'_\infty} \L_X^{-\ord (K_{Y'/Y} + g^*D)} d\m^X.
$$
\end{thm}

Thanks to the change of variables formula and Hironaka's
resolution of singularities,
one can reduce all computations to the case
in which $D$ is a simple normal crossing divisor
\cite[Theorem~0.2 and Notation~0.4]{K-M}.
Throughout this paper, we will use the following notation:
if $E_i$, with $i \in J$, are the irreducible components
of a simple normal crossing $\Q$-divisor on a nonsingular variety $Y$,
then for every subset $I \subseteq J$ we write
$$
E^0_I :=
\begin{cases}
Y \setminus E &\text{if $I = \emptyset$,} \\
(\cap_{i \in I}E_i) \setminus (\cup_{j \in J \setminus I}E_j)
&\text{otherwise.}
\end{cases}
$$
Now, consider a simple normal crossing
effective divisor $D = \sum_{i \in J} a_i E_i$
on a nonsingular $X$-variety $Y$ (here $E_i$ are the irreducible components of $D$).
Then a simple computation shows that
\begin{equation}\label{eq-snc-motivic}
\int_{Y_\infty} \L_X^{-\ord (D)} d\m^X = \sum_{I \subseteq J}
\frac{\{E^0_I\}}{\prod_{i \in I} \{\P_X^{a_i}\}}.
\end{equation}
The computations are carried out (over $\Spec \C$) in \cite{Cra}.
The importance of this formula is not just computational.
In fact, it implies that every integral of the form~(\ref{int})
is an element in the image of the natural ring homomorphism
$$
\r : K_0(\Sch_X)[\{\P_X^a\}^{-1}]_{a \in \N} \to \^\MM_X.
$$

\begin{rmk}
It is not known whether $\r$ is injective.
\end{rmk}

We let
\begin{equation}\label{N_X}
\NN_X := {\rm Im}(\r) \subset \^\MM_X.
\end{equation}
All together, we have a commutative diagram
$$
\xymatrix{
K_0(\Sch_X) \ar[d] \ar[dr]^\t \ar[r] & \MM_X \ar[d] \\
K_0(\Sch_X)[\{\P_X^a\}^{-1}]_{a \in \N} \ar[r]_(.7)\r &\^\MM_X.
}
$$
Note that the image of $\t$ is contained in $\NN_X$.

Given a morphism $h : V \to X$ of complex varieties, we obtain
a ring homomorphism $\ff_h : K_0(\Sch_X) \to K_0(\Sch_V)$
such that $\ff_h(\{Y\}) := \{Y \times_X V\}$ for every
$X$-scheme $Y$ (see \cite[Section~4]{Loo}).
We have $\ker(\t_X) \subseteq \ker(\t_V\o\ff_h)$,
hence a commutative diagram of ring homomorphisms
\begin{equation}\label{h}
\xymatrix{
K_0(\Sch_X) \ar[d]_{\t_X} \ar[r]^{\ff_h} &K_0(\Sch_V) \ar[d]^{\t_V} \\
\^\MM_X \ar[r] & \^\MM_V.
}
\end{equation}
If $h : V \to X$ is \'etale,
then we denote the image in $\^\MM_V$ of an element $\a \in \^\MM_X$ by $\a|_V$.
The following property will be used in the proof of Theorem~\ref{MK-Phi}.

\begin{prop}\label{restriction}
Let $h : V \to X$ be an \'etale morphism of complex varietes.
Let $D$ be an effective divisor on a nonsingular $X$-variety $Y$,
and let $D_V := D \times_X V$ and $Y_V := Y \times_X V$.
Then $Y_V$ is a nonsingular $V$-variety, $D_V$ is an effective divisor on $Y_V$, and
$$
\left(\int_{Y_\infty} \L_X^{-\ord (D)} d\m^X\right)\Big{|}_V
= \int_{(Y_V)_\infty} \L_V^{-\ord (D_V)} d\m^V.
$$
\end{prop}

\begin{proof}
Note that $Y_V$ is an open subset of $Y$ and $D_V = D|_V$.
Let $\p_m : Y_\infty \to Y_m$ and $\s_m : (Y_V)_\infty \to (Y_V)_m$
the truncation maps. Via the natural inclusions
$(Y_V)_\infty \subseteq Y_\infty$ and $(Y_V)_m \subseteq Y_m$,
we have $\s_m = \p_m|_{(Y_V)_m}$. Then
for any $p,m \in \N$ we have
\begin{align*}
\s_m(\Cont^p(D_V))
&= \s_m( \{\g \in (Y_V)_\infty \mid \ord_\g(D_V) = p \}) \\
&= \p_m(\{\g \in Y_\infty \mid \g(0) \in Y_V,\; \ord_\g(D) = p \} ) \\
&= \p_m(\Cont^p(D))\times_Y Y_V \\
&= \p_m(\Cont^p(D))\times_X V.
\end{align*}
Note also that $\L_V = \L_X|_V$. Then the assertion
follows by the definition of integral~\eqref{int}.
\end{proof}

\section{From the Relative Motivic Ring to Constructible Functions}

In this section we introduce a natural way
to read off constructible functions on a fixed variety $X$
out of motivic integrals that are computed relative to $X$.
A group homomorphism from relative Grothendieck rings to
groups of constructible functions also appears in \cite{BSY}, and
a general theory of motivic integration and constructible
functions has been announced by Cluckers and Loeser \cite{CL}.

Fix a complex algebraic variety $X$.
Let $F(X)$ be the group of constructible functions on $X$.
This is the subgroup of the abelian group of all $\Z$-valued functions
$f: X(\C)\rightarrow \Z$ that is generated by characteristic
functions $\1_S$, where $S$ ranges among closed subvarieties
of $X$. It is sometimes convenient to consider $F(X)$ as a ring,
with the product defined pointwise. For instance, for two constructible sets
$S$ and $T$ on $X$, we have $\1_S \. \1_T = \1_{S \cap T}$.
Then $F(X)$ is a commutative ring
with zero element $\1_\emptyset$ (the constant function ${\mathbf 0}$)
and identity $\1_X$ (the constant function ${\mathbf 1}$).
Basic properties of constructible functions may be found in \cite{Joy}.

Associated to any
morphism of varieties $f : Y \to X$, there is a group homomorphism
$f_* : F(Y) \to F(X)$ such that, for any
constructible set $S \subset Y$, the function $f_*\1_S$ is defined pointwise
by setting
$$
(f_*\1_S)(x) = \q_c(f^{-1}(x) \cap S) \quad\text{for every $x \in X$,}
$$
where $\q_c$ denotes the Euler characteristic
with compact support (in the analytic topology).
Since $f|_S$ is piecewise topologically locally trivial over a
stratification of $X$ in Zariski-locally closed subsets
(e.g., see \cite[Corollaire~(5.1)]{Ver}), the
function $f_*\1_S$ is indeed a constructible function on $X$.

Next recall the following well known properties.

\begin{lem}\label{GS}
\mbox{}
\begin{enumerate}
\item
If $Y = \bigsqcup Y_i$ is a decomposition of a complex
variety $Y$ as a disjoint union
of locally closed subvarieties, then $\q_c(Y) = \sum \q_c(Y_i)$.
\item
If $Z$ and $Z'$ are complex varieties, then
$\q_c(Z\times Z') = \q_c(Z)\q_c(Z')$.
\end{enumerate}
\end{lem}

It follows immediately by Lemma~\ref{GS}(a) that $f_*$ is a
group homomorphism \cite[Theorem~4.5]{Joy}.  We obtain a functor
$Y\mapsto F(Y), f\mapsto f_*$ from the category
of $X$-varieties to the category of abelian groups;
in particular, if $S$ is a subvariety of $Y$ and $g = f|_S$,
then $f_*\1_S = g_*\1_S$, where we view
$\1_S$ both as an element in $F(Y)$ and as an element in $F(S)$.

Note that, if $Y$ and $Y'$ are $X$-varieties and $x\in X$, then
$(Y\times_X Y')_x = Y_x \times Y'_x$. Thus, keeping in mind
Remark~\ref{rmk:sch-v-var}, it follows by Lemma~\ref{GS}(b) that
we can define a ring homomorphism $\Phi_0 : K_0(\Sch_X) \to F(X)$ by
setting
$$
\Phi_0(\{g:V \to X\}) = g_*\1_V
$$
for every $X$-variety $g : V \to X$,
and extending by linearity (see also \cite[Theorem~2.2]{BSY}). In the following,
recall the definition of $\NN_X$ from~\eqref{N_X}.

\begin{prop}\label{Phi}
There is a unique ring homomorphism $\Phi : \NN_X \to F(X)_\Q$ making
$$
\xymatrix{
K_0(\Sch_X) \ar[d]_\t \ar[r]^(.6){\Phi_0} & F(X) \ar@{^{(}->}[d] \\
\NN_X \ar[r]^\Phi & F(X)_\Q
}
$$
a commutative diagram of ring homomorphisms.
\end{prop}

\begin{proof}
Since $\Phi_0(\{\P_X^a\}) = (a+1)\1_X$ is an invertible element in
$F(X)_\Q$, $\Phi_0$ extends, uniquely, to a ring homomorphism
$$
\~\Phi :  K_0(\Sch_X)[\{\P_X^a\}^{-1}]_{a \in \N} \to F(X)_\Q.
$$
We claim that
\begin{equation}\label{ker}
\ker(\r) \subseteq \ker(\~\Phi).
\end{equation}
Let us grant this for now.
We conclude that $\~\Phi$ uniquely induces, and is uniquely
determined by, a ring homomorphism
$\Phi : \NN_X \to F(X)_\Q$. The commutativity of the
diagram in the statement is clear by the construction.

It remains to prove~(\ref{ker}). By clearing denominators
it is sufficient to show that, if $\a$ is
in the kernel of $\t$,
then $\Phi_0(\a) = 0$. When $X = \Spec \C$, this is proven in~\cite[(6.1)]{DL1}.
In general, consider an arbitrary morphism $x : \Spec \C \to X$.
By change of base, we obtain the commutative diagram~\eqref{h}
with $V = \Spec \C$. Since $\t(\a) = 0$,
this gives $\tau_{\Spec \C}(\psi_x(\a)) = 0$, hence $\chi_c(\psi_x(\a)) = 0$
by~\cite{DL1}. This means that $\Phi_0(\a)(x) = 0$.
Varying $x$ in $X$, we conclude that $\Phi_0(\a) = 0$.
\end{proof}

Given an effective divisor $D$ on a nonsingular $X$-variety $Y$, we define
$$
\Phi_{(Y,-D)}^X := \Phi \left(\int_{Y_\infty} \L_X^{-\ord (D)} d\m^X\right) \in F(X)_\Q.
$$
We will use the abbreviated notation $\Phi_{(Y,-D)}$ to denote this
function anytime $X$ is clear from the context. Moreover, if $D=0$, then
we write $\Phi_Y^X$, or just $\Phi_Y$.
We get at once the following properties:

\begin{prop}\label{restriction-Phi}
With the above notation, if $V$ is an open subset of $X$,
and we let $Y_V = Y \times_X V$ and $D_V = D|_{Y_V}$, then
$$
\Phi_{(Y_V,-D_V)}^V = \Phi_{(Y,-D)}^X \big{|}_V.
$$
\end{prop}

\begin{proof}
Apply $\Phi$ to both sides of the formula in Proposition~\ref{restriction}.
\end{proof}

\begin{prop}\label{change-of-vars-Phi}
Consider two nonsingular $X$-varieties $Y\xrightarrow{f}X$ and
$Y'\xrightarrow{f'}X$,
and assume that there is a proper birational morphism
$g : Y' \to Y$ over $X$.
Let $D$ be an effective divisor on $Y$.
Then
$$
\Phi_{(Y,-D)} = \Phi_{(Y',-(K_{Y'/Y}+g^*D))}
$$
in $F(X)_\Q$, where $K_{Y'/Y}$ is the relative canonical divisor of $g$.
\end{prop}

\begin{proof}
Apply $\Phi$ to both sides of the formula in
Theorem~\ref{change-of-vars}.
\end{proof}

\begin{cor}\label{cor-snc}
With the notation as in Proposition~\ref{change-of-vars-Phi},
assume that $K_{Y'/Y} + g^*D = \sum_{i \in J} a_iE_i$
is a simple normal crossing divisor. Then
\begin{equation}\label{snc}
\Phi_{(Y,-D)} = \sum_{I \subseteq J}
\frac{f'_*\1_{E^0_I}}{\prod_{i \in I} (a_i+1)}
\end{equation}
in $F(X)_\Q$.
In particular, if $D=0$, then $\Phi_Y = f_*\1_Y$.
\end{cor}

\begin{proof}
By Proposition~\ref{change-of-vars-Phi} and~\eqref{eq-snc-motivic}.
\end{proof}

\section{Constructible Functions Arising from Log-Terminal Pairs}

In the following we will use some terminology
coming from the theory of singularities of pairs;
standard references are \cite{Kol,K-M}.
As in the previous sections, we fix a complex variety $X$.
The goal of this section is to generalize the construction
introduced in the previous section and define a way to associate a
constructible function on $X$ to any {\it Kawamata log-terminal
pair} $(Y,\D)$, namely, a pair consisting of a normal $X$-variety $Y$
and a $\Q$-Weil divisor $\D$ on it such that $K_Y + \D$ is $\Q$-Cartier
and the pair has Kawamata log-terminal singularities.
We stress that $\D$ is not assumed here to be effective.

One needs to extend the motivic ring in order to
integrate order functions with rational values,
and so we now review this extension (see \cite{Vey1,Loo}).
We start considering the case in which
$Y$ nonsingular and $\D$ is a simple normal crossing $\Q$-divisor on $Y$.
We write
$$
D := - \D = \sum_{i \in J} a_i E_i,
$$
where $E_i$ are the irreducible components of $\D$.
In this case, the assumption of log-terminality for $(Y,\D)$
is equivalent to having $a_i > -1$ for all $i$.
We choose an integer $r$ such that $ra_i \in \Z$ for every $i$,
and define the ring $\^\MM_X^{1/r}$ to be the completion
of
$$
K_0(\Sch_X)[\L_X^{\pm 1/r}]
$$
with respect to a similar dimensional filtration as the one used in
the case $r=1$. Here $\L_X^{1/r}$ is a formal variable
with $(\L_X^{1/r})^r = \L_X$, and we assign to it dimension
$\frac 1r + \dim X$. Then we define
\begin{equation}\label{int-r}
\int_{Y_\infty} \L_X^{-\ord (D)} d\m^X
:= \sum_{p} \m^X(\Cont^p(rD))\.(\L_X^{1/r})^{-p}.
\end{equation}
(One can think that one is integrating
$(\L_X^{1/r})^{-\ord (rD)}$ instead of $\L_X^{-\ord (D)}$.)
Since $\Cont^p(rD)$ is non-empty only for integral values of $p$,
the summation appearing in the right hand side of~(\ref{int-r})
is taken over $\Z$. In fact, an explicit (and rather standard) computation
shows that the summation
is taken over $\N$ (this is not clear {\em a priori} because $D$ need not be effective).
This is the crucial point in order to ensure that
the integral defined above is indeed an element of $\^\MM_X^{1/r}$;
it is precisely at this point where we need the assumption
of log-terminality.
In addition, the same computation gives us the following
formula for the integral:
\begin{equation}\label{Q-gorenstein fmla}
\int_{Y_\infty} \L_X^{-\ord(D)} d\m^X = \sum_{I \subseteq J}\{E^0_I\}\prod_{i \in I}
\frac{\sum_{t=0}^{r-1} (\L_X^{1/r})^t}{\sum_{t=0}^{r(a_i+1)-1} (\L_X^{1/r})^t}.
\end{equation}
(One clearly sees from this formula that the assumption $a_i > -1$
is necessary for the summations in the denominators to be non-empty.)

The expression~\eqref{Q-gorenstein fmla}
is an element in the image of the ring homomorphism
$$
K_0(\Sch_X)\left[\left(\sum_{t=0}^{rd-1}
(\L_X^{1/r})^t\right)^{-1}\right]_{d \in \frac 1r\N^*} \longrightarrow \^\MM_X^{1/r}.
$$
We denote by $\NN_X^{1/r}$ the image of this homomorphism.
We extend the ring homomorphism $\Phi_0 : K_0(\Sch_X) \to F(X)$,
defined in the previous section, to a ring homomorphism
$$
\Phi_0 : K_0(\Sch_X)[\L_X^{1/r}] \to F(X)
$$
by setting $\Phi_0(\L_X^{1/r}) = \1_X$.
Observing that $\Phi_0(\sum_{t=0}^b (\L_X^{1/r})^t) = (b+1)\1_X$,
we conclude (as in the proof of Proposition~\ref{Phi}) that
$\Phi_0$ induces a ring homomorphism
$$
\Phi : \NN_X^{1/r} \to F(X)_\Q.
$$
Note that, for every rational number $a > -1$ and any choice of $r$
such that $ra \in \Z$,
$$
\Phi\left(\frac{\sum_{t=0}^{r-1}
(\L_X^{1/r})^t}{\sum_{t=0}^{r(a+1)-1} (\L_X^{1/r})^t}\right)
= \frac r{r(a + 1)}\1_X = \frac {\1_X}{a+1},
$$
which, in particular, does not depend on the choice of $r$.
Therefore, we can define
$$
\Phi_{(Y,\D)}^X := \Phi\left(\int_{Y_\infty} \L_X^{-\ord (D)} d\m^X\right) \in F(X)_\Q.
$$
By~\eqref{Q-gorenstein fmla} and the above discussion,
this function does not depend on the choice
of the integer $r$ needed to compute it.

\begin{rmk}\label{snc-lt}
The formula stated in~(\ref{snc})
still holds in the setting of this section, namely allowing $a_i$
to be rational numbers larger than $-1$.
\end{rmk}

Bearing in mind \cite[Lemma~2.30]{K-M},
Proposition~\ref{change-of-vars-Phi} extends to this setting,
giving us the following property.

\begin{prop}\label{change-of-vars-lt}
Consider a Kawamata log-terminal
pair $(Y,\D)$, with $Y$ a nonsingular $X$-variety.
Let $g : Y' \to Y$ be a proper birational morphism
such that $Y'$ is nonsingular and $\D' := - K_{Y'/Y} + g^*\D$
is a simple normal crossing
divisor. Then $(Y',\D'))$ is a log-terminal pair, and
$$
\Phi_{(Y,\D)}^X= \Phi_{(Y',\D')}^X.
$$
\end{prop}

We are now ready to consider the general setting: we start with a
Kawamata log-terminal pair $(Y,\D)$ over $X$,
namely, a pair consisting of a normal $X$-variety $Y$
and a $\Q$-Weil divisor $\D$ on $Y$ such that $K_Y + \D$ is $\Q$-Cartier
and the pair $(Y,\D)$ has Kawamata log-terminal
singularities \cite[Definition~2.34]{K-M}.
We can find a resolution of singularities $g : Y' \to Y$ such that, if
$\Ex(g)$ is the exceptional locus of $g$
and $\D' \subset Y'$ is the proper transform of $\D$, then
$\Ex(g) \cup \D'$ is a simple normal crossing divisor on $Y'$.
A resolution of this type is called a {\it log-resolution}.
Fix an integer $m$ such that $m(K_Y + \D)$ is Cartier.
Then there exists a unique $\Q$-divisor $\G$
on $Y'$ such that $m\G$ is linearly equivalent
to $- mK_{Y'} + g^*(m(K_Y + \D))$ and $\G + \D'$ is supported on $\Ex(g)$;
furthermore, $\G$ does not depend on the choices of
$K_Y$ and $m$ \cite[Section~2.3]{K-M}.
Note that $\G$ is a simple normal crossing $\Q$-divisor
by our assumption on the resolution,
and that the pair $(Y',\G)$ is Kawamata log-terminal \cite[Lemma~2.30]{K-M}).
Thus we can define
$$
\Phi_{(Y,\D)}^X := \Phi_{(Y',\G)}^X \in F(X)_\Q.
$$
By Proposition~\ref{change-of-vars-lt},
this definition is independent of the choice of resolution.
Similar abbreviations of the notation as in the previous section
will be used.

A particular case of this construction is given
by the following situation. Let $Y$ be a $X$-variety,
and suppose that $Y$ is normal, $K_Y$ is $\Q$-Cartier,
and the singularities of $Y$ are log-terminal
(that is, that $(Y,0)$ is a Kawamata log-terminal pair).
For short, we will say that $Y$ has (at most) {\it log-terminal singularities}.
Then we take a resolution of singularities $Y' \to Y$
with simple normal crossing exceptional divisor.
The relative canonical $\Q$-divisor $K_{Y'/Y}$ is then uniquely defined,
and we obtain the function
$$
\Phi_Y := \Phi_{(Y',-K_{Y'/Y})}^X \in F(X)_\Q
$$
by the above construction.

\begin{rmk}\label{rmk:beyond-lt:Phi}
While (\ref{Q-gorenstein fmla}) may not make sense is some $a_i$
is less than $-1$, the formal expression given in (\ref{snc}) does
make sense as soon as all $a_i \ne -1$; however, if some $a_i$ is
less than $-1$, then it is an open question whether the resulting
expression for $\Phi_{(Y,\D)}^X$ is, in general, independent of
the resolution (cf. \cite[Question~I]{Vey}). Nevertheless, if $Y$
is quasi-projective, $\D = 0$, and such condition on the
discrepancies $a_i$ is satisfied by a log-resolution of $Y$
factoring through the blowing up of the Jacobian ideal of $Y$,
then it follows by a result of Veys \cite{Vey3} that the formal
expression in (\ref{snc}) (with $D = 0$) does not depend on such
resolution. See also \cite{Vey2} for related results in dimension
two.
\end{rmk}

\section{Stringy Chern Classes}

In this section we define stringy Chern classes and
prove certain basic properties of these. Given a complex variety $X$,
MacPherson \cite{Mac} defined a homomorphism of additive groups
$$
c : F(X) \to A_*(X)
$$
such that, when $X$ is smooth, one has $c(\1_X) = c(T_X) \cap [X]$.
Abusing notation, we will also denote the extension
of this homomorphism to rational coefficients
$F(X)_\Q \to A_*(X)_\Q$ by $c$.
MacPherson also proved that these homomorphisms commute with
push-forwards along proper morphisms.

When $X$ is singular, MacPherson uses the above transformation to define a
generalization of total Chern class of $X$, by considering
$c(\1_X)$. The same class was independently defined by Schwartz \cite{Sch}.
This class is denoted by $c_{\SM} (X)$,
and is called the {\it Chern-Schwartz-MacPherson class} of $X$.
An easy consequence of functoriality of MacPherson transformation
is that, if $X$ is a proper variety, then $\int_X c_{\SM}(X) = \q(X)$.

We use the results of previous section to define
an alternative generalization of total Chern class of a singular variety,
more in the spirit of the newly introduced string-theoretic
invariants arising from motivic integration.
We assume that $X$ has at most log-terminal singularities
(in particular, $X$ is normal and $K_X$ is $\Q$-Cartier).
Via the construction given in the previous sections,
we obtain an element $\Phi_X \in F(X)_\Q$ associated to
the pair $(X,0)$.

\begin{defi}\label{definition}
Let $X$ be variety with at most log-terminal singularities.
The {\em stringy Chern class} of $X$ is the class
$$
c_{\st}(X) := c(\Phi_X) \in A_*(X)_\Q.
$$
If $\Phi_X$ is a integral-valued function, then we also
define the class
$$
\~c_{\st}(X) := c(\Phi_X) \in A_*(X).
$$
\end{defi}

\begin{rmk}
This definition can be given even if $X$ is a normal varieties
with $\Q$-Cartier canonical divisor and singularities worst
than log-terminal, if $X$ is quasi-projective and admits a log-resolution
factoring through the blowing up of the Jacobian ideal of $X$ such that
all discrepancies $a_i$ appearing in this resolution
are different from $-1$ (see Remark~\ref{rmk:beyond-lt:Phi}).
\end{rmk}

\begin{rmk}
If $\Phi_X \in F(X)$, then the image of $\~c_{\st}(X)$
in $A_*(X)_\Q$ is equal to $c_{\st}(X)$.
In general, let $X$ be a variety as in Definition~\ref{definition},
and let $m$ be a natural number such that,
for some resolution of singularities $Y \to X$ with
simple normal crossing relative canonical $\Q$-divisor
$K_{Y/X} = \sum k_i E_i$,
$m$ is divisible by every product $\prod_{i \in I} (k_i + 1)$
for which $E_I^0 \ne \emptyset$.
Then the stringy Chern class of $X$ can actually be defined
as an element in $A_*(X)_\LLL$, where $\LLL = \frac 1m\Z$.
\end{rmk}

In \cite{Bat2}, Batyrev defines the stringy Euler number $e_{\st}(X,\D)$
of a Kawamata log-terminal pair $(X,\D)$.
When $\D = 0$, this number is simply denoted by $e_{\st}(X)$.

\begin{prop}\label{e_st}
Let $X$ be a proper variety with at most
log-terminal singularities. Then
$$
\int_X c_{\st}(X) = e_{\st}(X).
$$
\end{prop}

\begin{proof}
Let $f : Y \to X$ be a resolution of singularities such that
$K_{Y/X} = \sum_i k_i E_i$ is a simple normal crossing $\Q$-divisor.
By Proposition~\ref{change-of-vars-Phi} and~(\ref{snc})
(see also Remark~\ref{snc-lt}) and the properness of $f$, we have
$$
c_{\st}(X) = c(\Phi_X) =  c(\Phi_{(Y,-K_{Y/X})}^X)
= c\left(\sum_{I \subseteq J} \frac{f_*\1_{E_I^0}}{\prod_{i \in I}(k_i+1)}\right)
= f_*c\left(\sum_{I \subseteq J} \frac{\1_{E_I^0}}{\prod_{i \in I}(k_i+1)}\right).
$$
Since the closure of each stratum of the stratification
$Y = \bigsqcup_{I \subseteq J} E^0_I$
is a union of strata, we can find rational numbers $b_I$
such that
$$
\sum_{I \subseteq J} \frac{\1_{E_I^0}}{\prod_{i \in I}(k_i+1)} =
\sum_{I \subseteq J} b_I \,\1_{\ov{E_I^0}}.
$$
Thus, keeping in mind that MacPherson transformation is
a group homomorphism and recalling Lemma~\ref{GS}(a), we obtain
$$
\int_X c_{\st}(X) = \int_Y \sum_{I \subseteq J} b_I \, c_{\SM}(\ov{E_I^0})
= \sum_{I \subseteq J} b_I  \int_{\ov{E_I^0}}  c_{\SM}(\ov{E_I^0})
= \sum_{I \subseteq J} b_I \, \q(\ov{E_I^0})
= \sum_{I \subseteq J} \frac{\q_c(E_I^0)}{\prod_{i \in I}(k_i+1)}.
$$
Now, we just observe that
the last side of this chain of equalities is precisely $e_{\st}(X)$
(see \cite[Definition~1.4]{Bat2}).
\end{proof}

When $X$ admits a crepant resolution (that is,
a resolution of singularities $Y \to X$ with trivial relative
canonical divisor), the stringy Chern class
of $X$ has been already defined in \cite{Alu} (although it was not named this
way). The equivalence of the definitions follows by the next statement.

\begin{prop}\label{crepant}
If $X$ admits a crepant resolution $f : Y \to X$,
then $c_{\st}(X) = f_*(c(T_Y) \cap [Y])$. In particular,
if $X$ is smooth, then $c_{\st}(X) = c(T_X) \cap [X]$.
\end{prop}

\begin{proof}
Since $K_{Y/X} = 0$, $f$ is proper, and $Y$ is nonsingular, we have
$$
c_{\st}(X) = c(\Phi_X) = c(\Phi_Y^X) = c(f_*\1_Y) = f_*c(\1_Y)
= f_*(c(T_Y) \cap [Y]).
$$
\end{proof}

We extend Kawamata's definition
of $K$-equivalence \cite[Definition 1.1]{Kaw} as follows.

\begin{defi}\label{defi-K}
Two normal varieties $X$ and $X'$ with $\Q$-Cartier canonical divisors
are said to be {\it $K$-equivalent}
if there exists a nonsingular variety $Y$ and proper and birational morphisms
$f : Y \to X$ and $f' : Y \to X'$ such that $K_{Y/X} = K_{Y/X'}$
(as divisors).
\end{defi}

\begin{rmk}\label{resolutions}
If $X$ and $X'$ are $K$-equivalent varieties, then
the condition on the relative canonical divisors
is satisfied for every choice of $Y$, $f$ and $f'$.
\end{rmk}

\begin{rmk}
If $X$ and $X'$ have at most canonical singularities, then
the condition $K_{Y/X} = K_{Y/X'}$ in Definition~\ref{defi-K}
can be replaced with the condition $K_{Y/X} \equiv K_{Y/X'}$.
Indeed in this case the two conditions are equivalent by the
Negativity Lemma \cite[Lemma~3.39]{K-M}.
In particular, if both $X$ and $X'$ are smooth,
then both conditions are equivalent to
$f^*\O_X(K_X) \cong {f'}^*\O_{X'}(K_{X'})$.
\end{rmk}

Using the formalism of motivic integration,
we generalize the main result of \cite{Alu}.

\begin{thm}\label{K}
Let $X$ and $X'$ be varieties with at most log-terminal singularities,
and assume that they are $K$-equivalent. Consider any diagram
$$
\xymatrix{
&Y \ar[dl]_f \ar[dr]^{f'} \\
X & &X',
}
$$
with $Y$ a nonsingular variety and $f$ and $f'$ proper and birational morphisms.
Then:
\begin{enumerate}
\item  There is a class $C \in A_*(Y)_\Q$ such that
$f_*C = c_{\st}(X)$ in $A_*(X)_\Q$ and $f'_*C = c_{\st}(X')$ in $A_*(X')_\Q$.
\item
Moreover $K_{Y/X} =K_{Y/X'}$ and,
assuming that this is equal to a simple normal crossing
divisor $\sum_{i \in J} k_iE_i$
(here the $E_i$ are the irreducible components),
then $C$ is the class
$$
C = c\left(\sum_{I \subseteq J} \frac{\1_{E_I^0}}{\prod_{i \in I}(k_i+1)}\right).
$$
\end{enumerate}
\end{thm}

\begin{proof}
By assumption, $K_{Y/X} = K_{Y/X'}$. Let $K$ denote this divisor.
It is enough to prove the theorem assuming that $K$ has simple normal
crossings. Indeed, by further blowing up $Y$, we can always reduce to this
case, and push-forward on Chow rings is functorial for
proper morphisms. Then, defining $C$ as in part~(b) of the statement, we have
$$
f_*C = c\left(\sum_{I \subseteq J} \frac{f_*\1_{E_I^0}}{\prod_{i \in I}(k_i+1)}\right)
= c(\Phi_{(Y,-K)}^X) = c(\Phi_X) = c_{\st}(X),
$$
where we have applied the functoriality of $c$ with respect
to proper morphisms for the first equality,
used~(\ref{snc}) for the second one, and applied
Proposition~\ref{change-of-vars-Phi} for the third.
Similarly, $f'_* C = c_{\st}(X')$.
\end{proof}

The case when $X$ and $X'$ are smooth has already been considered in~\cite{Alu}.

\begin{cor}[Aluffi]
Let $X$ and $X'$ be smooth $K$-equivalent varieties.
Let $Y$ be a nonsingular variety with proper birational
morphisms $f : Y \to X$ and $f': Y \to X'$ as in Theorem~\ref{K}.
Then there is a class $C \in A_*(Y)_\Q$ such that
$f_*C = c(T_X) \cap [X]$ in $A_*(X)_\Q$
and $f'_*C = c(T_{X'}) \cap [X']$ in $A_*(X')_\Q$.
\end{cor}

\begin{proof}
By Theorem~\ref{K} and Proposition~\ref{crepant}.
\end{proof}

It is possible to carry out a general theory
of stringy Chern classes for pairs by defining
$$
c_{\st}(X,\D) := c(\Phi_{(X,\D)}) \in A_*(X)_\Q
$$
for any Kawamata log-terminal pair $(X,\D)$.
Then an obvious adaptation of the proof of Proposition~\ref{e_st} gives
$$
\int_X c_{\st}(X,\D) = e_{\st}(X,\D).
$$

\begin{rmk}\label{rmk-BL}
If $Y$ is a smooth variety and $- \D = \sum_{i\in J} a_i E_i$
is a simple normal crossing divisor on $Y$ with rational
coefficients $a_i > -1$, then the adjunction formula \cite[Example~3.2.12]{Ful},
together with an inclusion-exclusion argument, gives
$$
c_{\st}(Y,\D) = (c(T_Y) \cap [Y]) \.
\prod_{i \in J} \frac{1 + \frac{1}{a_i+1}[E_i]}{1 + [E_i]}
$$
in $A_*(Y)_\Q$.
It was explained to us by Borisov and Libgober
how the expression in the right hand side can be obtained
from the orbifold elliptic class, defined in
\cite[Definition~3.2]{BL4},
by taking the coefficients of a certain
Laurent expansion in one of the two variables after
a limiting process in the other variable.
In particular, the stringy Chern class of a variety with at most
log-terminal singularities can
also be reconstructed from its orbifold elliptic class.
\end{rmk}

\section{The Grothendieck Ring of Deligne-Mumford Stacks}\label{stacks}

In this section we freely use the language of stacks. General references
for the theory of stacks include \cite{LMB} and \cite{Gomez}.
The goal in this section is to define the Grothendieck
ring of Deligne-Mumford stacks (in the relative setting),
and to prove certain fundamental properties of this ring
that will be useful in the study of stringy Chern classes of quotient varieties.
We will work in greater generality than is
actually needed in the remainder of this paper: it seems useful
 not only to give more conceptual proofs of these
properties, but also to develop the more general framework for possible
future use.
The construction of the Grothendieck ring of stacks also appears
in a new version of the preprint \cite{Joy2}.
See also \cite{Yas1,Yas2} for motivic integration over stacks.

Fix a base variety $X$ over $\C$. We denote by $\DM_X$ the category of
separated Deligne-Mumford stacks of finite type over $X$,
hereafter {\it DM stacks over $X$}.
We also let
$\Space_X$ denote the category of separated algebraic spaces of finite
type over $X$.
Let $K_0(\DM_X)$ (respectively, $K_0(\Space_X)$)
 be the quotient of the free abelian group generated
by isomorphism classes of objects of $\DM_X$ (objects of $\Space_X$)
 over $X$, modulo the relations
$\{\XX \setminus \YY \} = \{ \XX \} - \{ \YY \}$ whenever $\YY \to \XX$
is a closed immersion. Then $K_0(\DM_X)$ and $K_0(\Space_X)$ become
commutative rings with identity under the product defined by
fiber product over $X$.

Recall the following definition.

\begin{defi}
Let $\XX$ be a DM stack over $X$. The {\it inertia stack}
$I \XX$ of $\XX$ is defined as follows.
An object of $I \XX$ is a pair $(\x, \alpha)$
with $\x$ an object of $\XX$ and $\alpha \in \Aut(\x)$ (note that
$\a$ is an automorphism of $\x$ over $X$). A morphism
$(\x, \alpha) \to (\z, \beta)$ in $I \XX$ is a morphism
$\phi \colon \x \to \z$ in $\XX$ (hence over $X$) with $\alpha \phi = \phi \beta$.
\end{defi}
The following lemma gives an alternative characterization of the inertia stack
which will be used in the proof of Proposition \ref{inertia prop}.
\begin{lem}\label{lemma-IX}
For a DM stack $\XX$ over $X$, we have
\begin{equation}\label{IX}
I\XX = \XX \times_{\SS}\XX,
\end{equation}
where $\SS = \XX\times_X \XX$ and both maps
$\XX\rightarrow \SS$ are the diagonal morphism.
\end{lem}

\begin{proof}
By \cite[Lemma~1.13]{Vistoli}, the inertia stack $I\XX$, viewed as a stack
over $\C$ via the forgetful functor $\DM_X \to \DM_\C$, is equal to
$$
\XX \times_{\XX \times_{\C} \XX}\XX,
$$
where both maps $\XX\rightarrow \XX\times_{\C}\XX$ are the diagonal morphism.
We check that
$$
\XX \times_{\XX \times_{\C} \XX}\XX =
\XX \times_{\XX\times_X \XX}\XX
$$
by applying the universal property in both directions.
By the natural morphism $\XX\times_X\XX \to \XX\times_{\C} \XX$,
we get a morphism from the right hand side to the left one.
The inverse morphism is constructed similarly by using the
natural morphism from the diagonal of $\XX\times_\C\XX$ to $\XX\times_X\XX$.
\end{proof}

The inertia stack comes equipped with a representable morphism
$I\XX\rightarrow \XX$.

\begin{prop}\label{inertia prop}
The inertia map $\XX \mapsto I\XX$ induces a ring homomorphism
\begin{displaymath}
I \colon K_0(\DM_X) \to K_0(\DM_X).
\end{displaymath}
\end{prop}

\begin{proof}
It suffices to check that:
\begin{enumerate}
\item
If $\YY\subset \XX$ is a closed substack, then the natural morphism
$I\YY \rightarrow I\XX$ is a closed immersion, and
$I(\XX\setminus\YY) = (\XX\setminus\YY)\times_\XX I\XX =
I\XX\setminus I\YY$.
\item
If $\XX$ and $\XX'$ are DM stacks over $X$, then
$I(\XX\times_X\XX') = I\XX \times_X I\XX'$.
\end{enumerate}
Note that, by definition \cite[Remarque 3.5.1]{LMB}, a closed substack
$\YY\subset\XX$
is a full subcategory with the property that if $\x$ is an object of $\YY$
and $\z$ is an object of $\XX$ isomorphic to $\x$ then $\z$ is also an object
of $\YY$.
It then follows immediately from the definitions
that $I\YY =  I\XX \times_\XX \YY$, and, since
closed immersions pull back to closed immersions,
the first part of (a) follows.
The statement about $I(\XX\setminus\YY)$ follows from the fact that
the morphism $I\YY\rightarrow I\XX\times_\XX \YY$ is a bijection.

Part (b) is immediate from \eqref{IX} of Lemma \ref{lemma-IX}: we have
\begin{align*}
I(\XX\times_X\XX')
&= (\XX\times_X\XX')\times_{(\XX\times_X\XX')\times_X(\XX\times_X\XX')}
(\XX\times_X\XX') \\
&= (\XX\times_{\XX\times_X\XX}\XX)\times_X(\XX'\times_{\XX'\times_X\XX'}\XX')
= I\XX\times_X I\XX'.
\end{align*}
\end{proof}

There is a natural ring homomorphism
$$
\Psi : K_0(\Sch_X) \to K_0(\Space_X)
$$
given by the obvious inclusion of categories.  Moreover, because
every algebraic space has an open set that is a scheme, $\Psi$ is
actually an isomorphism. Indeed, the inverse of $\Psi$ is defined
by taking, for every element in $K_0(\Space_X)$, a representative
in the free $\Z$-module generated by isomorphism classes of
objects in $\Space_X$, and then stratifying into locally closed
subschemes every algebraic space whose isomorphism class appears
in this formal combination. Common refinement of two such
stratifications shows that the corresponding element in
$K_0(\Sch_X)$ is independent of the particular stratification
chosen.

\begin{prop}\label{coarse ring map}
The functor that to each object $\XX$ of $\DM_X$ associates its coarse
moduli space $\coarse{\XX}$ defines,
after composition with $\Psi^{-1}$, a ring homomorphism
\begin{displaymath}
K_0(\DM_X)\longrightarrow K_0(\Sch_X).
\end{displaymath}
\end{prop}

\begin{proof}
To prove the assertion, there are two properties to verify:
\begin{enumerate}
\item[(a)]
If $\YY\subset \XX$ is a closed immersion of objects of $\DM_X$ with open
complement $\UU = \XX \setminus \YY$, then
$[\coarse{\YY}] + [\coarse{\UU}]= [\coarse{\XX}]$ in $K_0(\Space_X)$.
\item[(b)]
If $\XX_1$ and $\XX_2$ are objects of $\DM_X$, then the coarse moduli spaces
satisfy $[\coarse{\XX_1\times_X\XX_2}] = [\coarse{\XX}_1\times_X\coarse{\XX}_2]$.
\end{enumerate}
We will rely on the following lemma.
\begin{lem}\label{bijections give eq}
Let $Z\xrightarrow{f} Y$ be a morphism in $\Space_X$.  Suppose that $f$ induces
a bijection $Z(\xi)\rightarrow Y(\xi)$ for every geometric point $\xi$.  Then
$[Z]=[Y]$ in $K_0(\Space_X)$.
\end{lem}
\begin{proof}[Proof of Lemma]
Suppose first that $Z$ and $Y$ are schemes.  After stratifying as
necessary, we obtain a collection $\{Z_i \rightarrow Y_i\}$ such that
each $f_i$ is flat and bijective on geometric points, in particular a
homeomorphism.  A further stratification and replacing schemes by their
associated reduced schemes then guarantees that each $f_i$ is an \'etale
homeomorphism, hence an isomorphism.  Thus, in $K_0(\Sch_X)$ we get
$[Z]  = \sum [Z_i] = \sum [Y_i] = [Y]$.

Now, suppose that $Z$ and $Y$ are algebraic spaces, and stratify as above
so that each $f_i$ is flat.  We may
then stratify the $Z_i$ and $Y_i$ further so that each stratum
 is an object of $\Sch_X$; the flatness of each
$f_i$ guarantees that $f_i$ is open, moreover, and so these
further refinements may be made compatibly with the morphisms
$f_i$.  We thus obtain stratifications $Z=\coprod Z_i$, $Y=\coprod
Y_i$ into objects of $\Sch_X$ such that $f$ restricts to morphisms
$f_i: Z_i \rightarrow Y_i$ which are bijections on geometric
points. The conclusion of the last paragraph then guarantees the
conclusion of the lemma.
\end{proof}
We now return to the proof of the proposition.  By Corollary 1.3 and
Definition 1.8(G) of \cite{KeelMori}, the natural maps
$\XX\rightarrow \coarse{\XX}$ and $\UU\rightarrow \coarse{\UU}$ induce
bijections on geometric points.  It follows that the natural map
$\coarse{\YY}_{\operatorname{red}}\rightarrow
\coarse{\XX}\setminus\coarse{\UU}$ induces bijections on geometric points, and
thus Lemma \ref{bijections give eq} proves (a).

For (b), there is again a natural morphism
$\coarse{\XX_1\times_X \XX_2}\rightarrow \coarse{\XX}_1\times_X\coarse{\XX}_2$,
and it suffices to prove that this induces isomorphisms on geometric points.
To do this, observe that
$\coarse{\XX_1\times_X \XX_2}(\xi) = (\XX_1\times_X\XX_2)(\xi)$ by
\cite{KeelMori}, and that the latter is just
$\XX_1(\xi)\times_{X(\xi)}\XX_2(\xi)$, which also agrees with
$\coarse{\XX}_1(\xi)\times_{X(\xi)}\coarse{\XX}_2(\xi)
=(\coarse{\XX}_1\times_X\coarse{\XX}_2)(\xi)$ since the natural maps
$\XX_i(\xi)\rightarrow \coarse{\XX}_i(\xi)$ are bijections.
\end{proof}

In the remainder of this section, we focus on the case of orbifolds.
For a finite group $G$, we denote by $\CC(G)$ a set of representatives
of the conjugacy classes of elements in $G$, and by
$\SS(G)$ a set of representative of the conjugacy classes
of subgroups of $G$; the conjugacy class of an element $g$ of $G$ will
be denoted by $(g)$.
We denote the centralizer in $G$ of an element $g \in G$  by $C_G(g)$, or simply
by $C(g)$ when $G$ is clear from the context.
For a subgroup $H$ of $G$,
$N_H$ will denote the normalizer of $H$ in $G$.
We remark that every representative of a
conjugacy class in $N_H$ of an element of $H$
is also in $H$. In particular, there is no ambiguity in the notation
$\CC(N_H) \cap H$.

\begin{lem}\label{inertia}
For a quotient stack $[M/G]$ over $X$, we have
$$
I [M/G] \cong \bigsqcup_{g \in \CC(G)} [M^g/C(g)],
$$
where $M^g \subset M$ is the fixed point set of $g$.
\end{lem}

\begin{proof}
By definition
$$
I[M/G] = \left[\left( \bigsqcup_{g \in G} M^g\times\{g\} \right)\Big{/}G \right] =
\bigsqcup_{g \in \CC(G)} \left[\left(\bigsqcup_{g' \in (g)} M^{g'}
\times \{g'\} \right) \Big{/}G \right]
$$
where $h \in G$ acts on $M^g \times \{g\}$ by taking
$(x ,g )$ to $(xh, h^{-1}gh) \in M^{ h^{-1}gh} \times \{ h^{-1}gh\}$.
Now one can see that
$$
\left[\left(\bigsqcup_{g' \in (g)} M^{g'} \times \{g'\} \right)\Big{/}G \right]
\cong [M^g/C(g)].
$$
\end{proof}

For any subgroup $H$ of $G$, set
$$
M^H = \{ y \in M \mid G_y = H\}.
$$
Here $G_y$ is the stabilizer of $y$ in $G$.
Note that $M^H$ is contained in the fixed locus of $H$, but may be
strictly smaller. As $H$ runs among all subgroups of $M$,
we get a stratification of $M$, that we can write as follows:
$$
M = \bigsqcup_{H \in \SS(G)} \left(\bigsqcup_{H' \in (H)} M^{H'}\right).
$$

\begin{prop}\label{identity-stacks}
The following equality holds in $K_0(\DM_X)$:
$$
\{ I[M/G] \} =
\sum_{H \in \SS(G)} \left(\sum_{h\in\CC(N_H)\cap H} \{[M^H/C_{N_H}(h)]\}\right).
$$
\end{prop}

\begin{proof}
We start observing that
$$
\{[M/G]\} = \sum_{H \in \SS(G)}
\left\{\left[\left(\bigsqcup_{H' \in (H)} M^{H'}\right)\Big{/}G \right] \right\}.
$$
For $g \in G$ and $x \in M^{H'}$, $xg$ belongs to $M^{g^{-1} H' g}$, so
$$
\left[\left(\bigsqcup_{H' \in (H)} M^{H'}\right)\Big{/}G \right] \cong
[M^H/N_H].
$$
Then the proposition follows by  Lemma~\ref{inertia},
after we apply the inertia map
to the stacks $[M^H/N_H]$ and observe that, for a given $g \in N_H$,
$(M^H)^g$ is equal to $M^H$ if $g \in H$, and empty otherwise.
\end{proof}

We deduce the following fact from
Propositions~\ref{identity-stacks} and \ref{coarse ring map}.

\begin{cor}\label{identity-McKay}
Let $M$ be a variety with an action of a finite group $G$,
and let $X = M/G$.
Then the following identity holds in $K_0(\Sch_X)$:
$$
\sum_{g \in \CC(G)}\{M^g/C(g)\} =
\sum_{H \in \SS(G)} \left(\sum_{h\in\CC(N_H)\cap H} \{M^H/C_{N_H}(h) \}\right).
$$
\end{cor}

\section{McKay Correspondence for Stringy Chern Classes}\label{section-MK}

In this section we compare stringy Chern classes of quotient varieties
with Chern-Schwartz-MacPherson classes of fixed-point set data.

Let $M$ be a smooth quasi-projective complex variety of dimension $d$,
and let
$G$ be a finite group with an action on $M$.
Let $X = M/G$, with projection $\p : M \to X$.
We will assume that
the canonical line bundle $\omega_M$
of $M$ descends to the quotient: that is, we assume that for every
$p\in M$, the natural
action of the stabilizer $G_p$ on $\bigwedge^d T_p^*(M)$
is trivial.
This is for instance a natural setup considered in mirror symmetry;
a motivating example is that of a finite subgroup of $\SL_n(\C)$
acting on $\C^n$.
(The assumption that $M$ is quasi-projective can be dropped
by simply requiring that the action of $G$ is {\it good},
as defined in \cite[Section~5]{Loo}.)
Under these hypotheses, $X$ is normal
(e.g., see \cite[Proposition~2.3.11]{DK})
and has Gorenstein canonical singularities
(e.g., see \cite[Subsection~1.3]{Rei}); it also follows that
$\p^*\omega_X = \omega_M$.

There are two ways of ``breaking the orbifold into simpler pieces''.
The first way is to stratify $X$ according
to the stabilizers of the points on $M$.
For any subgroup $H$ of $G$, let
$X^H \subseteq X$ be the set of points $x$
such that, for every $y \in \p^{-1}(x)$,
the stabilizer of $y$ is conjugate to $H$.
If we let $H$ run in a set $\SS(G)$ of representatives of
conjugacy classes of subgroups of $G$, we obtain a
stratification of $X$:
$$
X = \bigsqcup_{H \in \SS(G)} X^H.
$$
The second way is to look at the fixed-points sets
as orbifolds under the action of the corresponding
centralizers.
For every $g \in G$, consider the fixed-point set $M^g \subseteq M$---note
that it is smooth since $M$ is smooth.
There is a commutative diagram
$$
\xymatrix{
M^g \ar[d] \ar@{^{(}->}[r] & M \ar[d]^\p \\
M^g/C(g) \ar[r]^(.6){\p_g} & X.
}
$$
It is easy to see that $\{M^g/C(g)\}$, as an element
in $K_0(\Sch_X)$, is independent of the representative $g$
chosen for its conjugacy class in $G$.
Note also that $\p_g$ is a proper morphism.
In the following, we fix a set $\CC(G)$ of representatives
of conjugacy classes of elements of $G$.

\begin{thm}\label{MK-Phi}
With the above assumptions and notation, $\Phi_X$ is an element in $F(X)$
and the following identities hold in $F(X)$:
\begin{equation}\label{eq-MK-Phi}
\Phi_X = \sum_{H\in\SS(G)}|\CC(H)|\.\1_{X^H}
= \sum_{g \in \CC(G)} (\p_g)_*\1_{M^g/C(g)}.
\end{equation}
\end{thm}

\begin{proof}
Let $f : Y \to X$ be any resolution of singularities,
and fix a common multiple $r$ of the denominators of the
coefficients of $K_{Y/X}$.
Note that for every open subset $V \subseteq X$, the restriction
of $f$ to $Y_V := f^{-1}(V)$ is a resolution of singularities of $V$,
and that $K_{Y/X}|_{Y_V} = K_{Y_V/V}$. Then, by Proposition~\ref{restriction-Phi},
it suffices to check~\eqref{eq-MK-Phi} on each piece of an
open covering of $X$. Therefore we can assume without loss
of generality that $\om_M = \O_M$ and that there exists a non-vanishing
$G$-invariant global section of $\om_M$.
Then the ``motivic McKay correspondence" \cite{Bat2,DL2,Rei,Loo},
in the relative setting, gives the following identity in $\NN_X$:
\begin{equation}\label{X^H}
\int_{Y_\infty} \L_X^{-\ord (K_{Y/X})} d\m^X
= \sum_{H\in\SS(G)}\{X^H\}\.\(\sum_{h\in\CC(H)}\L_X^{\age(h)}\).
\end{equation}
(See in particular the formulation given in \cite[Theorem~4.4(4)]{Rei}.)
Here $\age(h)$ are non-negative integers depending on the action
of $h$ on $M$ (these numbers are called {\it ages}, see \cite{Rei}).
Actually, one may need to split the first sum in~\eqref{X^H}
further to separate different connected components of $X^H$,
as the ages might change from one component to the other,
but we will apply $\Phi$ to~\eqref{X^H} in a moment,
and these numbers therefore will play no role.

Applying $\Phi : \NN_X \to F(X)_\Q$ to both sides of~(\ref{X^H}),
we get the first identity stated in the theorem.
In particular, this shows that $\Phi_X \in F(X)$.
In order to prove the second identity, we apply $\Phi$ to both sides
of the equation in Corollary~\ref{identity-McKay}. This gives
$$
\sum_{g \in \CC(G)} (\p_g)_* \1_{M^g/C(g)} =
\sum_{H \in \SS(G)} \left(\sum_{h\in\CC(N_H)\cap H} \Phi(\{M^H/C_{N_H}(h) \})\right).
$$
Fix $H \in \SS(G)$ and $h \in H$.
Note that $C_H(h) = C_{N_H}(h) \cap H$, and consider the commutative
diagram with exact rows
$$
\xymatrix{
1 \ar[r] &C_H(h) \ar@{^{(}->}[d] \ar[r]
&C_{N_H}(h) \ar@{^{(}->}[d] \ar[r] &K_H(h) \ar@{^{(}->}[d] \ar[r] &1 \\
1 \ar[r] &H \ar[r] &N_H \ar[r] &K_H \ar[r] &1.
}
$$
Since $H$ acts trivially on $M^H$, whereas $K_H$ acts freely,
we have $M^H/N_H = M^H/K_H$, and the quotient map $M^H \to M^H/K_H$
is \'etale. Similarly, $M^H/C_{N_H}(h) = M^H/K_H(h)$.
Thus we have a commutative diagram
\begin{equation}\label{1}
\xymatrix{
M^H \ar[d]_{\p|_{M^H}} \ar[r]^(.4)\e
&M^H/C_{N_H}(h) \ar@{=}[r] \ar[d]^\n & M^H/K_H(h) \ar[d] \\
X^H \ar@{=}[r] &M^H/N_H \ar@{=}[r] &M^H/K_H
}
\end{equation}
where all maps are \'etale. Moreover,
\begin{equation}\label{2}
\deg \n = \frac{\deg (\p|_{M^H})}{\deg \e} = \frac{|K_H|}{|K_H(h)|}.
\end{equation}
Now, let $(h)_{N_H}$ be the conjugacy class of $h$ in $N_H$,
and $(h)_H$ be the class of $h$ in $H$.
Then, to conclude, it is enough to observe that
\begin{equation}\label{3}
\begin{split}
\sum_{h\in\CC(N_H)\cap H} \frac{|K_H|}{|K_H(h)|}
&= \sum_{h\in\CC(N_H)\cap H}
\frac{|N_H|}{|C_{N_H}(h)|}\Big{/}\frac{|H|}{|C_H(h)|}\\
&= \sum_{h\in\CC(N_H)\cap H}\frac{|(h)_{N_H}|}{|(h)_H|}\\
&= \sum_{h \in H} \frac 1{|(h)_H|} = |\CC(H)|.
\end{split}
\end{equation}
Combining~\eqref{1},~\eqref{2} and~\eqref{3}, we conclude that
$$
\sum_{h\in\CC(N_H)\cap H} \Phi(\{M^H/C_{N_H}(h) \}) = |\CC(H)|\.\1_{X^H}
$$
in $F(X)$.
Thus, taking the sum over $\SS(G)$, we obtain the
second identity stated in the theorem.
\end{proof}

The first identity in Theorem~\ref{MK-Phi} can be interpreted
in terms of the {\it orbifold push-forward}
defined by Joyce for constructible functions on stacks \cite{Joy}.
As this is quite unrelated to the rest
of the paper, we will not introduce the terminology
of constructible functions and push-forwards on stacks,
for which we refer to \cite{Joy}.

\begin{cor}\label{MK-orb}
With the notation introduced at the beginning
of this section, $\Phi_X$
is equal to the Joyce's orbifold push-forward of the constructible
function $\1_{[M/G](\C)}$ on $[M/G]$
along the morphism of stacks $[M/G] \to X$.
\end{cor}

We omit the proof of this corollary.
The reader familiar with Joyce's paper will see that the
proof follows directly from the various definitions
\cite[Definitions~5.3,~5.12 and~5.16]{Joy}.

Now we use the second identity in Theorem~\ref{MK-Phi} to
compare the stringy Chern class of $X$ with the Chern-Schwartz-MacPherson classes
of the fixed-point sets.
We remark that in general the quotients $M^g/C(g)$ may have several
irreducible components, but the construction of MacPherson
extends to this case, hence there is no problem in defining
the Chern-Schwartz-MacPherson class of $M^g/C(g)$.

\begin{thm}\label{MK-chern}
With the notation as in the beginning of this section, we have
$$
\~c_{\st}(X) = \sum_{g \in \CC(G)} (\p_g)_*c_{\SM}(M^g/C(g))
$$
in $A_*(X)$.
\end{thm}

\begin{proof}
Applying $c : F(X) \to A_*(X)$ to the first and last members
of the formula in Theorem~\ref{MK-Phi}, we obtain
$$
\~c_{\st}(X) = \sum_{g \in \CC(G)} c((\p_g)_*\1_{M^g/C(g)}).
$$
Hence the statement follows by recalling that $c$ commutes with $(\p_g)_*$.
\end{proof}

In their influential papers \cite{DHVW},
Dixon, Harvey, Vafa and Witten define the {\it orbifold Euler number}
of a quotient variety. Keeping the notation introduced in this
section and using the formulation of Hirzebruch and H\"ofer \cite{HH},
this number was defined as
$$
e(M,G) := \sum_{g \in \CC(G)} e_c(M^g/C(g)).
$$
where $e_c$ stands for Euler characteristic with compact supports.
It was observed in \cite{DHVW} that in some cases when $X$ admits
a crepant resolution $Y \to X$, this number is equal to the Euler
number of $Y$. This was then proved and generalized by Batyrev
(see \cite[Theorem~7.5]{Bat2} for the general result).

\begin{thm}[Batyrev]
With the notation as in the beginning of this section,
$$
e_{\st}(X) = e(M,G).
$$
\end{thm}

One can also see this formula from Theorem~\ref{MK-Phi}
by considering the constant morphism
$g : X \to \Spec\C$ and applying $g_*$ to the first and last terms of the
formula \eqref{eq-MK-Phi}. Indeed applying $g_*$ to $\Phi_X$ we recover
$e_{\st}(X)$, whereas applying it to the last term we obtain $e(M,G)$.
When $X$ is proper, this can also be viewed
from Theorem~\ref{MK-chern} by taking degrees and applying
Proposition~\ref{e_st}.

\providecommand{\bysame}{\leavevmode \hbox \o3em
{\hrulefill}\thinspace}

\end{document}